\documentclass{amsart} 
\usepackage[all]{xy}
\numberwithin{equation}{section}
\title[Explicit models]{Explicit models for some stable
  categories of maximal Cohen-Macaulay modules} 
\author{Louis de
  Thanhoffer de V\"olcsey} 
\author{Michel Van den Bergh}
\usepackage{amssymb}
\let\cal\mathcal

\def\Cscr{{\cal C}}

\def\Pscr{{\cal P}}

\let\blb\mathbb

\def \ZZ{{\blb Z}}

\def \NN{{\blb N}}

\def\Lotimes{\overset{L}{\otimes}}

\def\mod{\operatorname{mod}}

\def\rad{\operatorname {rad}}

\def\Spec{\operatorname {Spec}}

\def\Hom{\operatorname {Hom}}

\def\End{\operatorname {End}}
\def\RHom{\operatorname {RHom}}

\def\coker{\operatorname {coker}}

\def\Tor{\operatorname {Tor}}
\def\End{\operatorname {End}}

\def\r{\rightarrow}

\def\MCM{\operatorname{MCM}}
\def\cone{\operatorname{cone}}
\def\Perf{\operatorname{Perf}}

\newtheorem{lemma}{Lemma}[section]
\newtheorem{proposition}[lemma]{Proposition}
\newtheorem{theorem}[lemma]{Theorem}

\newtheorem{lemmas}{Lemma}[subsection]
\newtheorem{propositions}[lemmas]{Proposition}
\newtheorem{theorems}[lemmas]{Theorem}

\theoremstyle{definition}

\newtheorem{examples}[lemmas]{Example}

{

}

\theoremstyle{remark}

\newtheorem{remarks}[lemmas]{Remark}

\newdimen\uboxsep \uboxsep=1ex
\def\uboxn#1{\vtop to 0pt{\hrule height 0pt depth 0pt\vskip\uboxsep
\hbox to 0pt{\hss #1\hss}\vss}}

\def\uboxs#1{\vbox to 0pt{\vss\hbox to 0pt{\hss #1\hss}
\vskip\uboxsep\hrule height 0pt depth 0pt}}

\email{louicious@mac.com, michel.vandenbergh@uhasselt.be}
\address{Departement WNI,  Universiteit Hasselt,
 3590 Diepenbeek, Belgium.}
\thanks{The second author is a senior researcher at the FWO}
\begin{document}
\begin{abstract}
We give concrete DG-descriptions of certain stable categories
of maximal Cohen-Macaulay modules. This makes in possible to describe the latter as generalized cluster categories in some cases.
\end{abstract}

\maketitle
\section{Introduction}
Throughout $k$ is a field.  Let $(R,m)$ be a local complete noetherian
Gorenstein $k$-algebra of Krull dimension $n$ with an isolated
singularity and with $R/m\cong k$. Then it is well-known that the
stable category of Cohen-Macaulay $R$-modules $\underline{\MCM}(R)$ is
an $n-1$-Calabi-Yau category. This category is also sometimes called
the singularity category of $R$.

Assume given a Cohen-Macaulay $R$-module $N$ such that
$\Lambda=\End_R(R\oplus N)$ has finite global dimension and is a Cohen-Macaulay $R$-module. Thus
$\Lambda$ is a non-commutative crepant resolution of $R$ in the sense
of \cite{VdB32}. Under these conditions Iyama has shown that $N$ is an
$n-1$-cluster tilting object in $\underline{\MCM}(R)$
\cite[Thm 5.2.1]{iyama} and it is a natural question if $\underline{\MCM}(R)$ may
be obtained by one of the standard constructions of CY-categories with
cluster tilting object \cite{Amiot,Guo}.

\medskip

Assume for simplicity that $N$ contains no repeated
summands and no summands isomorphic to $R$. Let $(S_i)_{i=0,\ldots,l}$
be the simple $\Lambda$-modules with~$S_0$ corresponding to the summand
$R$ of $R\oplus N$. Let $e_0$ the idempotent given by the
projection $R\oplus N\r R$. Put $l=\Lambda/\rad\Lambda$ and $l^0=le_0$.

In this note we observe the following (Lemma \ref{ref-4.1-3}, Theorem \ref{ref-5.1.1-5}, Remark
\ref{ref-5.1.2-6})
\begin{theorem}
\label{ref-1.1-0}
Let $(T_lV,d)\r \Lambda$ be a finite minimal model for $\Lambda$ (see
\S\ref{secpreliminaries}).
%
Put $\Gamma=T_lV/T_lV e_0T_lV$. Then one has
\[
\underline{\MCM}(R)\cong \Perf(\Gamma)/\langle (S_i)_{i\neq 0}\rangle
\]
and furthermore $\Gamma$ has the following properties:
\begin{enumerate}
\item $\Gamma$ has finite dimensional cohomology in each degree;
\item $\Gamma$ is 
concentrated in degrees $\le 0$;
\item $H^0(\Gamma)=\Lambda/\Lambda e_0\Lambda$;
\item
As a graded algebra $\Gamma$ is of the form
$T_{l^0} V^0$ for $V^0=(1-e_0)V(1-e_0)$.
\end{enumerate}
\end{theorem}
The miminal model $(T_lV,d)$ always exists for abstract reasons (see Lemma \ref{ref-4.1-3} below)
but to use Theorem~\ref{ref-1.1-0} effectively one must be able to
describe it. There are at least two cases where this is easy:
\begin{itemize}
\item $\Lambda$ is $3$-Calabi-Yau. It then follows from \cite{VdBsuper} that the minimal
model of $\Lambda$ is derived from a super potential.
\item $\Lambda$ is Koszul.
\end{itemize}

We obtain in particular (see Proposition \ref{ref-5.2.1-16} below)
\begin{proposition} \
\label{ref-1.2-1} Assume that $\Lambda$ is 3-Calabi Yau. 
Then 
$\underline{\MCM}(R)$ is the generalized cluster category associated
to a quiver with super potential as introduced by Claire Amiot in \cite{Amiot}. 
\end{proposition}
The Koszul case is well taylored to quotient singularities. 
In \S\ref{ref-6.4-23} we consider the case where
$R=k[[x_1,\ldots,x_n]]^G$ for $G\cong \ZZ/m\ZZ$ acting linearly
(Proposition \ref{ref-6.5.1-26}) with weights $(\xi^{a_i})_i$
for $\xi$ a primitive $m$'th root of unity and $0\le a_i\le m-1$, $\gcd(a_i,m)=1$, $\sum_i a_i\cong 0$ mod $m$. We
given an explicit description of $\Gamma$ in this case.

As an application (Proposition \ref{ref-6.6.1-28} below) of this
explicit description of $\Gamma$ we obtain a new proof of the following result from
\cite{AIR,iyama2} 
\begin{proposition}
\label{ref-1.3-2}
Assume $\sum_i a_i=m$.
  Let $P$ be the quiver with vertices $[m-1]\overset{\text{def}}{=}\{1,\ldots,m-1\}$ and for
  each $j\in [m-1]$ and $i\in [n]$ such that $j+a_i\in [m-1]$ an arrow
  $x_i:j\r j+a_i$ with relations $x_ix_{i'}=x_{i'}x_i$. Let $C$ be the
  path algebra of~$P$ modulo these relations. Let
  $\Theta_C=\RHom_{C^e}(C,C\otimes C)$ be the inverse dualizing
  complex of~$C$ \cite{Keller11,VdB29}. Then there is a quasi-isomorphism of DG-algebras
\[
\Gamma\cong T_C(\Theta_{C}[n-1])
\]
In particular $
\underline{\MCM}(R)
$
is a generalized $n-1$-cluster category.
\end{proposition}
\section{Acknowledgement}
As said Proposition \ref{ref-1.3-2} was proved independently in
\cite{AIR} (see also \cite{iyama2}).  Some time ago the authors had been looking at explicit
DG-models for certain stable categories of maximal Cohen-Macaulay
models (based on Theorem \ref{ref-1.1-0}). We knew that Proposition
\ref{ref-1.3-2} was true in the 3-dimensional case (as then it is
a consequence of Proposition \ref{ref-1.2-1}) and expected a
suitable generalization to hold in higher dimension as well. However
we became distracted by other interests.

In May 2010 the authors learned about the precise statement of Proposition
\ref{ref-1.3-2} from a very interesting lecture by Osamu Iyama at the OberWolfach
meeting ``Interactions between Algebraic Geometry and Noncommutative
Algebra'' (see \cite{iyama2}). In that lecture the result was derived from
a very general statement about graded Calabi-Yau algebras of Gorenstein parameter 1.

We expected the result to be a consequence of Theorem
\ref{ref-1.1-0} as well and this is indeed the case although the 
derivation is a little more involved than we anticipated. Nonetheless since the methods
in \cite{AIR} appear to be quite different from ours we decided it would
be useful to write down our own proof.

We are very grateful to Dong Yang for his careful reading of the first
version of this manuscript. His comments have helped us to improve the exposition.
Dong Yang has also informed us that he and Martin Kalck have a different approach to some of
our results based on ``recollement''.
\section{Notation and conventions}
Throughout $l$ is a separable algebra over the ground field $k$. If we
refer to an ``algebra'' $\Gamma$ then we mean an $l$-algebra. In other
words $\Gamma$ is implicitly equipped with a ring homomorphism $l\r
\Gamma$. Unadorned tensor products are over $k$.
\section{Preliminaries}
\label{secpreliminaries}
In this paper we deal with two settings. The ``complete'' case and the
``graded'' case. This mainly affects our definition of a tensor
algebra. Let $V$ be an $l$-bimodule and denote by $T^\circ_lV$ the
ordinary tensor algebra of $V$. In the complete case $T_lV$ is the
completion of $T_l^\circ V$ at the twosided ideal $(T^\circ_lV)_+$
generated by $V$. If $V$ carries an additional $\ZZ$-grading (in a homological sense) then $T_l V$ refers to the graded
completion of $T^\circ_lV$ at $(T^\circ_lV)_+$. In the graded case $V$
will be equipped with an \emph{additional} (non-homological) grading which we sometimes
refer to as the ``Adams grading'' (following \cite{LPWZ}).  In this case
$T_lV$ is simply $T_l^\circ V$ equipped with the extended Adams
grading.
\begin{itemize}
\item In the complete case we deal with (DG-)algebras which can be
  written as of $T_lV/I$ with $V$ a finitely  generated (graded)
  $l$-bimodule and $I$ a (graded) closed ideal in $T_lV$.  
\item In the graded case we deal with graded (DG-)algebras which can be
  written as of $T_lV/I$ with $V$ a finitely generated
  (bi)graded $l$-bimodule and $I$ a (bi)graded ideal. 
Unless otherwise specified the Adams grading will be left bounded in each
homological degree. 
\end{itemize}
Below we will use finite minimal models for some
algebras.  For us a finite minimal model of an algebra $\Lambda$
will be a DG-algebra of the form $(T_lV,d)$ for a finitely generated
graded $l$-bimodule $V$ living in degrees $\le 0$ such that $dV$ is in the
twosided ideal generated by $V\otimes_l V$, together with a
quasi-isomorphism $T_lV\r\Lambda$ (with $\Lambda$ being viewed as a DG-algebra
concentrated in degree zero).
\begin{lemma} 
\label{ref-4.1-3} Let $\Lambda$ be either complete of the from $l+\rad \Lambda$ or graded
of the form $\Lambda=l+\Lambda_1+\cdots$. Assume $\dim \sum_i \Tor^\Lambda_i(l,l)<\infty$. Then
$\Lambda$ has a finite minimal model with 
\begin{equation}
\label{ref-4.1-4}
V=\bigoplus_{i\ge 1} \Tor^\Lambda_i(l,l)[i-1]
\end{equation}
Conversely if $\Lambda$ has a finite minimal model $T_lV$ then $V$ is given by the
formula \eqref{ref-4.1-4}.
\end{lemma}
\begin{proof} This follows from the bar-cobar formalism. In the graded
  case it has been proved explicitly in \cite{LPWZ} (see \S1.1 in
  loc.\ cit.). One checks that the proof also goes through in the
  complete case (see e.g.\ \cite[Prop. A.5.4]{VdBsuper}).
\end{proof}
For our examples we only need the case where $\Lambda$ is a Koszul algebra
or a completion thereof. In that case it is trivial to construct the minimal
model directly (see Proposition \ref{ref-6.1.1-18} and Remark \ref{ref-6.1.2-20} below). 
\section{Stable categories of Cohen-Macaulay modules}
In this section we show how to obtain explicit models for some
stable categories of Cohen-Macaulay modules. In order
to have somewhat concise statements we intentionally do not
work in the greatest possible generality. 
\subsection{General results}
In this section we assume  that $(R,m)$ is a complete local noetherian Gorenstein
$k$-algebra with an isolated singularity (by complete 
  local we mean either literally complete local or else $\NN$-graded local and $R/m=k$). Let
  $(M_i)_{i=1,\ldots,l}$ be indecomposable finitely generated maximal
Cohen-Macaulay 
  $R$-modules which are pairwise non-isomorphic and
also not isomorphic to $R$. Put $M_0=R$,
  $M=\bigoplus_i M_i$ and $\Lambda=\End_R(M)$. Assume that $\Lambda$ has
finite global dimension so that $\Lambda$ is a non-commutative crepant
resolution of $R$ in the sense of~\cite{VdB32}.

Let $e_i\in \Lambda$ be
  the idempotent given by the projection $M\r M_i$ and put $P_i=\Lambda e_i$.
 Let $S_i=P_i/\rad P_i$ be the
corresponding simple $\Lambda$-module. Put $l=\sum_i ke_i$. 
\begin{theorems}
\label{ref-5.1.1-5}
 Let $T_lV$ be a finite minimal model of $\Lambda$ (cfr Lemma \ref{ref-4.1-3}). Let $\underline{\MCM}(R)$ be the stable
  category of maximal Cohen-Macaulay $R$-modules (graded if we are in
  the graded context).  Put $\Gamma=T_lV/T_lV e_0 T_lV$. 
\begin{enumerate}
\item
 The $(S_i)_{i\neq 0}$ are perfect
  $\Gamma$-modules and furthermore there is an exact
  equivalence
\[
 \underline{\MCM}(R)=\Perf(\Gamma)/\langle (S_i)_{i\neq 0}\rangle
\]
\item The DG-algebra $\Gamma$ has finite dimensional
  cohomology in each degree. Furthermore $H^0(\Gamma)=\Lambda/\Lambda
  e_0\Lambda$.
\end{enumerate}
\end{theorems}
\begin{remarks} \label{ref-5.1.2-6} Put $V^0=(1-e_0)V(1-e_0)$ and
  $l^0=l/le_0$. Then it is not hard to see that
\[
T_lV/T_lV e_0 T_lV\cong T_{l^0}V^0
\]
\end{remarks}
\begin{proof}[Proof of Theorem \ref{ref-5.1.1-5}] We first prove (1). 
According to Lemma \ref{ref-5.1.3-8}  below we have
 an equivalence 
\[
D^b_f(\Lambda)/\langle (S_i)_{i\neq 0}\rangle\cong D^b_f(R):
\]
which sends $P_0$ to $R$.

Hence using Buchweitz's theorem \cite{BuchweitzAlone} we get
\begin{multline}
\label{ref-5.1-7}
\underline{\MCM}(R)=D^b_f(R)/\langle R\rangle=\biggl (D^b_f(\Lambda)/\langle (S_i)_{i\neq 0}
\biggr)/
\langle P_0\rangle\\
=D^b_f(\Lambda)/\langle P_0,(S_i)_{i\neq 0}\rangle
=\biggl(D^b_f(\Lambda)/\langle P_0\rangle\biggr)/\langle (S_i)_{i\neq 0}\rangle
\end{multline}
To obtain our conclusion it is now sufficient to prove 
\[
D^b_f(\Lambda)/\langle P_0\rangle=\Perf(T_lV/T_lV e_0 T_lV)
\]
By hypotheses $\Lambda$ has finite global dimension and hence $D^b_f(\Lambda)=\Perf(\Lambda)
=\Perf(T_lV)$.
We may now conclude as in \cite[Lemma 7.2]{Keller11} (and proof).

\medskip

Since $H^i(\Gamma)=\Hom_{\Perf(\Gamma)}(\Gamma,\Gamma[i])$ and 
$\Perf(\Gamma)\cong D^b_f(\Lambda)$ the second statement is an immediate consequence
of Proposition \ref{ref-5.1.4-11} below. 
\end{proof}
\begin{lemmas} 
\label{ref-5.1.3-8} 
The functor
\begin{equation}
\label{ref-5.2-9}
\Xi:D^b_f(\Lambda)\r D^b_f(R):N\mapsto e_0N
\end{equation}
induces an equivalence 
\[
D^b_f(\Lambda)/\langle (S_i)_{i\neq 0}\rangle\cong D^b_f(R):
\]
\end{lemmas}
\begin{proof}
  Let $U\in D^b_f(R)$. Since $P_0$ is a vector bundle on $\Spec R-\{m\}$ 
we have that $P_0\otimes_R U$ has finite dimensional cohomology.

Let $N$ be such that for $n\ge N$ we have $H^{-n}(U)=0$.  We claim
that for $n\ge N$ we have that $H^{-n}(P_0\Lotimes_R U)$ is an
extension of $(S_i)_{i\neq 0}$, i.e. $e_0 H^{-n}(P_0\Lotimes_R U)=0$.
Indeed
\begin{equation}
\label{ref-5.3-10}
e_0 H^{-n}(P_0\Lotimes_R U)=H^{-n}(e_0\Lambda e_0\Lotimes_R U)\cong
H^{-n}(R\Lotimes_R U)=H^{-n}(U)=0
\end{equation}
  Define $\Phi(U)=\tau_{\ge -N} ( P_0\Lotimes_R U)$. 
Then $\overline{\Phi(U)}$ is a well defined 
object of $D^b_f(\Lambda)/\langle (S_i)_{i\neq 0}\rangle$. We claim that
$\overline{\Phi(-)}$ yields a quasi-inverse to \eqref{ref-5.2-9}.

If $U\in D^b_f(R)$ then the computation \eqref{ref-5.3-10} shows that 
$\Xi\overline{\Phi}(U)=U$. Conversely assume $V\in D^b_f(\Lambda)$. Then
$\overline{\Phi}\Xi(V)=\tau_{\ge -N}(P_0\Lotimes_R e_0 V)$. Let $C$ be the cone
of $P_0\Lotimes_R e_0 V=\Lambda e_0\Lotimes_R e_0 V\r V$.  We see $e_0C=0$. 
In other words the cohomology of $C$ is given by extensions of $(S_i)_{i\neq 0}$.
Furthermore by our choice of $N$ we have $e_0H^{-n}(V)=H^{-n}(e_0V)=0$ 
for $n\ge N$ and hence $H^{-n}(V)$ is an also an extension 
of $(S_i)_{i\neq 0}$ for such $n$.
Thus working modulo $\langle (S_i)_{i\neq 0}\rangle$ we have
\[
\tau_{\ge -N}(P_0\Lotimes_R e_0 V)=\tau_{\ge -N}V=V
\]
which finishes the proof. 
\end{proof}
\begin{propositions}
\label{ref-5.1.4-11}
The category
$D^b_f(\Lambda)/\langle P_0\rangle$ is $\Hom$-finite and in addition
\begin{equation}
\label{ref-5.4-12}
\Hom_{D^b_f(\Lambda)/\langle P_0\rangle}(\Lambda,\Lambda)=\Lambda/\Lambda e_0\Lambda
\end{equation}
\end{propositions}
\begin{proof}
Since $\Lambda$ has finite global dimension to prove $\Hom$-finiteness
it is sufficient to prove for any $M\in D^b_f(\Lambda)$:
\begin{equation}
\label{ref-5.5-13}
\dim \Hom_{D^b_f(\Lambda)/\langle P_0\rangle}(\Lambda,M)<\infty
\end{equation}
If $N\in\mod(\Lambda)$ then there is a map
\begin{equation}
\label{phidef}
\phi:P_0^n\r N
\end{equation}
such that $\coker\phi\in \langle (S_i)_{i\neq 0}\rangle$.
Using $\phi$ we can ``resolve'' any $M\in D^b_f(\Lambda)$ with a complex
$P\in \langle P_0\rangle$ such that $\cone(P\r M)$ modulo $\langle
(S_i)_{i\neq 0}\rangle$ is an object $M_1$ in $\mod(\Lambda)[n]$ for
$n\gg 0$.  
Hence to prove
\eqref{ref-5.5-13} we have
reduced ourselves to the cases $M\in \mod(\Lambda)[n]$ for $n\gg 0$ or $M=S_i[p]$ for
$i\neq 0$ and $p\in \ZZ$.

To deal with the case $M=S_i[p]$ let $N$ be an extension of
$(S_i)_{i\neq 0}$ in $\mod(A)$. Then $\Hom^i_{D^b_f(\Lambda)}(P_0,N)=0$
from which we easily deduce
\begin{equation}
\label{ref-5.6-14}
\Hom_{D^b_f(\Lambda)/\langle P_0\rangle}^i(\Lambda,N)
=
\begin{cases}
N &\text{if $i=0$}\\
0& \text{otherwise}
\end{cases}
\end{equation}
Thus this case is clear. 

To deal with the other case it is sufficient to show the following statement.
Assume $M\in \mod(\Lambda)$ then 
\begin{equation}
\label{ref-5.7-15}
\Hom^i_{D^b_f(\Lambda)/\langle P_0\rangle}(\Lambda,M)=0
\end{equation}
for $i> 0$. 
Let $p$ be a map $\Lambda\r M[i]$ in $D^b_f(\Lambda)/\langle P_0\rangle$ with $i>0$. Then 
$p$ is represented in $D^b_f(\Lambda)$ by a diagram of the following kind
\[
\xymatrix{%
& C\ar[dr]^{p'}\ar[dl]_{q} &\\
\Lambda && M[i]
}
\]
such that $P\overset{\text{def}}{=}\cone q\in \langle P_0\rangle$. It
is easy to see that we may assume that this diagram is an actual
diagram of complexes and that in addition $C$ is a finite complex of
finitely generated projectives. Then it follows that we may also
assume that $P$ is a finite complex of finite direct sums of $P_0$.

The
composition $P[-1]\r C\xrightarrow{p'} M[i]$ factors through
$\sigma_{\le -i}(P[-1])$ where $\sigma_{\le -i}$ denotes ``naive'' truncation.
Hence $\sigma_{\le -i}(P[-1])\in \langle P_0\rangle$.

 We then obtain a morphism 
of distinguished triangles in $D^b_f(\Lambda)$
\[
\xymatrix{
P[-1]\ar[r]\ar[d] & C\ar[d]^{p'}\ar[r]^q\ar@{-}[d] & \Lambda\ar@{.>}[d]\\
\sigma_{\le -i}(P[-1])\ar[r]& M[i]\ar[r] & Z
}
\]
where $Z$ is the cone of the lower leftmost map. 
We see that $Z$ has no cohomology in degree zero and hence
$\Hom_{D^b_f(\Lambda)}(\Lambda,Z)=0$. Thus $p'$ factors
through $\sigma_{\le i}(P[-1])$.  So $p$ (which is the image of $p'$ in
$D^b_f(\Lambda)/\langle P_0\rangle$) is the zero map. This finishes
the proof of \eqref{ref-5.7-15} for $i>0$. 

It follows also from  \eqref{ref-5.7-15} that $\Hom_{D^b(\Lambda)/\langle
  P_0\rangle}(\Lambda,-)$ is a right exact functor on
$\mod(\Lambda)$.  To compute $\Hom_{D^b(\Lambda)/\langle
  P_0\rangle}(\Lambda,M)$ let $\bar{M}=\coker \phi$ (cfr \eqref{phidef}). Applying
\eqref{ref-5.7-15} together with \eqref{ref-5.6-14} to the right exact
sequence
\[
  P_0^n \r M \r \bar{M}\r 0
\]
we obtain 
\[
\Hom_{D^b(\Lambda)/\langle P_0\rangle}(\Lambda,M)=\bar{M}
\]
Applying this identity with
$M=\Lambda$ we immediately obtain \eqref{ref-5.4-12}.
\def\qed{}\end{proof}
\subsection{Application to Claire Amiot's generalized cluster categories}
Let $(Q,w)=(Q_0,Q_1,w)$ be a quiver with super potential (thus $w$ is
an element of the completed path algebra $\widehat{kQ}$). The
Ginzburg algebra $\Gamma(Q,w)$ is the DG-algebra
$(k\widetilde{Q},d)$ where $\widetilde{Q}$ is the graded quiver
with vertices $Q_0$ and arrows
\begin{itemize}
\item The original arrows $a$ in $Q_1$ (degree 0);
\item Opposite arrows $a^\ast$ for $a\in Q_1$ (degree -1);
\item Loops $c_i$ at vertices $i\in Q
 _0$ (degree -2). We put
  $c=\sum_i c_i$.
\end{itemize}
The differential is 
\begin{align*}
da&=0\qquad a\in Q_1\\
da^\ast&=\frac{\partial w}{\partial a}\qquad a\in Q_1\\
dc&=\sum_{a\in Q_1} [a^\ast,a]
\end{align*}
Let $\widehat{\Gamma}(Q,w)$ be the graded completion of $\Gamma(Q,w)$
at path length.  The cohomology in degree zero of $\widehat{\Gamma}(Q,w)$ is
called the Jacobian algebra of $(Q,w)$ and is denoted by $\Pscr(Q,w)$. Concretely
one has
\[
\Pscr(Q,w)=\widehat{kQ}\left/\left\langle \left(\frac{\partial w}{\partial
    a}\right)_{a\in Q_1}\right\rangle\right.
\]
Put
\[
\Cscr_{Q,w}=\Perf(\widehat{\Gamma}(Q,w))/\langle (S_i)_{i\in Q_0} \rangle
\] 
where $S_i$ is the simple
$k\widetilde{Q}$ representation corresponding to the vertex $i$ (this is
perfect because of \cite[\S2.14]{KY}).  In
\cite[Thm 3.6]{Amiot} Claire Amiot shows that if $\Pscr(Q,w)$ is finite
dimensional then $\Cscr_{Q,w}$ is a 2-CY category
with cluster tilting object $\widehat{\Gamma}(Q,w)$.
The endomorphism ring of this cluster tilting object is  $\Pscr(Q,w)$.

\medskip

Assume now that $R$ is a complete local 3-dimensional Gorenstein domain with 
residue field $k$ with an isolated singularity.
Assume that $N$
is a Cohen-Macaulay $R$-module whose summands are not isomorphic to
$R$ and pairwise non-isomorphic such that in addition the following property
holds:

\begin{itemize}
\item For $\Lambda=\End_R(R\oplus N)$ one has $\Lambda\cong H^\ast(\widehat{\Gamma}(Q,w))$
  for some finite quiver and potential $w$.
\end{itemize}
This condition is equivalent to
the 3-Calabi-Yau property for $\Lambda$ (see
\cite{Gi,Keller11,VdBsuper}).

Let $0$ be the vertex of $Q$ corresponding to the $R$-summand of $R\oplus N$. 
Let $(Q^0,w^0)$ be obtained from $(Q,w)$
by deleting all edges in $Q$ which are adjacent to $0$ and all paths
 in $w$ which pass through $0$. 
By Theorem \ref{ref-5.1.1-5}
one obtains
\begin{propositions}
\label{ref-5.2.1-16}
One has an exact equivalence of triangulated categories
\[
\underline{\MCM}(R)\cong \Cscr_{Q^0,w^0}
\]
\end{propositions}
\begin{examples}
Assume that $R=k[[t,x,y,z]]/(xy-zt)$. Then $R$ has a non-commutative crepant resolution
$\End_R(R\oplus I)$ with $I=(x,z)$. The corresponding quiver is
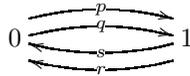
\[
\xymatrix{
0 
\ar@<1.5ex>@/^0.5em/[rr]|p\ar@/^0.5em/[rr]|q&&
1
\ar@<1.5ex>@/^0.5em/[ll]|r\ar@/^0.5em/[ll]|s
}
\]
with super potential
\[
w=psqr-prqs
\]
We see that $Q^0$ consists of the single vertex $1$ and $w^0=0$. Hence
\[
\widehat{\Gamma}(Q^0,w^0)=k[c]
\]
with $\deg c=-2$. 
Using the definitions in \cite{Keller6} it
follows that $\underline{\MCM}(R)$ is the cluster category associated
to the single vertex/no loops quiver; which is simply the category of $\ZZ/2\ZZ$
graded finite dimensional vector spaces. 
\end{examples}
\section{Minimal models for Koszul algebras}
The first few sections contain no new material. Their main purpose is to make explicit some formulas.
Throughout we follow a version of the Sweedler convention. An element
$a$ of a tensor product $X\otimes Y$ is written as $a'\otimes a''$ (thus we suppress
the summation sign) and a similar convention for longer tensor
products.
\subsection{The general case}
\label{ref-6.1-17}
Below $A=T_l^\circ V/(R)$ is a finitely generated quadratic algebra
(thus $V$ is a finitely generated $l$-bimodule and $R\subset V\otimes_l V$).  We put
\[
J_n=\bigcap_{i=0}^{n-1} V^{\otimes i-1}\otimes R\otimes V^{\otimes n-1-i} 
\]
and in particular
$J_1=V$, $J_2=R$.

We note that $A$ is naturally $\NN$-graded by giving $V$ degree $1$ (as usual this is referred to
as the Adams grading). 

The canonical map 
\[
J_n\mapsto J_i\otimes_l J_{n-i}
\]
is denoted by $\delta_i$ (or $\delta_{i,n-i}$ for clarity). It has degree zero for the Adams grading. 
Following our convention we write. 
\[
\delta_i(a)=\delta_i(a)'\otimes \delta_i(a)''
\]

Put 
\[
\tilde{V}=\bigoplus_{n>0} J_n[n-1]
\]
and
\[
\tilde{A}=T_l\tilde{V}
\]
For $a\in J_n$ considered as an element of degree $n-1$ in $\tilde{A}$
we put
\[
da=(-1)^{i-1}\sum_i \delta_i(a)'\otimes \delta_i(a)''
\]
It is easy to see that $d^2=0$.
\begin{propositions} 
\label{ref-6.1.1-18} Assume that $A$ is Koszul. Then the map $q:\tilde{A}\r A$  induced
from the projection map $\tilde{V}\r J_1=V$
is a quasi-isomorphism.
\end{propositions}
\begin{proof} This follows in the standard way from the bar-cobar
  formalism (see Lemma \ref{ref-4.1-4}). For the benefit of the reader we give a direct proof.

By definition the fact that $A$ is Koszul means that the following complex of graded left $A$-modules is exact
\begin{equation}
\label{ref-6.1-19}
\cdots \r A\otimes_l J_2\r A\otimes_l J_1\r A\r l\r 0
\end{equation}
where the differential is given by
\[
d:A\otimes_l J_n\r A\otimes_l J_{n-1}:a\otimes b\mapsto a\delta_{1,n-1}(b)'\otimes \delta_{1,n-1}(b)''
\]
Put
\[
M=(T_l\tilde{V})_+\overset{\text{def}}{=}\bigoplus_{n\ge 1} \tilde{V}^{\otimes_l
  n}\subset T_l\tilde{V}
\]
We consider $M$ as a left $\tilde{A}$ sub DG-module of $\tilde{A}$. As left graded $\tilde{A}$-module
we have
\[
M=\tilde{A}\otimes_l (J_1\oplus J_2[1]\oplus J_3[2]\oplus\cdots)
\]
Let $\tilde{C}$ be the cone of
the inclusion map $i:M\r \tilde{A}$. As graded $\tilde{A}$-modules we have
\[
\tilde{C}=\tilde{A}\otimes_l (l\oplus J_1[1]\oplus J_2[2]\oplus\cdots)
\]
As $\coker i=l$ the obvious map $\tilde{C}\r l$ is a a quasi-isomorphism.

Put $C=A\otimes_{\tilde{A}}\tilde{C}$. Then one checks that $C$ is
precisely the complex \eqref{ref-6.1-19} (without the right most $l$). Thus
$C\r l$ is a quasi-isomorphism as well and hence so is the canonical
map $\tilde{C}\r C$.

We now equip $\tilde{C}$ with an ascending filtration of
sub-DG-$\tilde{A}$-modules as follows: $F_0 \tilde{C}=\tilde{A}$,
$F_1\tilde{C}=\tilde{A}\otimes(k\oplus J_1[1])$, \dots. We equip $C$ with the
similar filtration.
The canonical
map $\tilde{C}\r C$ is a map of filtered DG-$\tilde{A}$-modules.

Assume we have shown that $\tilde{A}\r A$ is a quasi-isomorphism in
Adams degree $\le n$. Then $(\tilde{C}/F_0\tilde{C})_{n+1}\r (C/F_0C)_{n+1}$ is 
a quasi-isomorphism.
 Given that $\tilde{C}_{n+1}\r C_{n+1}$ 
is a quasi-isomorphism we deduce that $F_0\tilde{C}_{n+1}=\tilde{A}_{n+1}\r A_{n+1}=(F_0C)_{n+1}$ is also a
quasi-isomorphism.
\end{proof}
\begin{remarks} 
\label{ref-6.1.2-20} It is easy to see that a suitable analogue of Proposition \ref{ref-6.1.1-18}
holds for the completed rings $(T^\circ_lV/(R))\,\hat{}$. Indeed this simply means
that the directs sums involved in the grading become direct products. Taking
direct products of vector spaces is an extremely well behaved functor so
it does not break anything. We will rely on this fact below without further
elaboration. 
\end{remarks}
\subsection{Polynomial rings}
We recall the familiar finite minimal model for polynomial rings.
Assume $A=k[x_1,\ldots,x_n]$.  If we put $V=\sum_i kx_i$ then $J_n=\wedge^n V\subset V^{\otimes n}$.
Write $x_{S}=\wedge_{i\in S} x_i$ for $S\subset [n]$, $[n]=\{1,\ldots,n\}$, $S\neq \emptyset$. In the wedge product $\wedge_{i\in S} x_i$
we assume that the indices of the variables are in ascending order. 

Specializing \S\ref{ref-6.1-17} to this situation we obtain that
$\tilde{A}$ is equal to $k\langle(x_S)_{S\neq \emptyset}\rangle$ with differential
\begin{equation}
\label{ref-6.2-21}
dx_S=\sum_{S=A\coprod B,A\neq\emptyset,B\neq\emptyset}(-1)^{|A|-1} \epsilon_{A,B} x_A x_{B}
\end{equation}
where $\epsilon_{A,B}$ is the sign defined by
\[
\wedge_{i\in S} x_i=\epsilon_{A,B}\left(\wedge_{i\in A} x_i\right)\wedge \left(\wedge_{i\in B} x_i\right)
\]
\subsection{Crossed products}
We now assume that $k$ has characteristic zero.
Assume that $\Lambda=A\# G$ is a crossed product with $A=k[x_1,\ldots,x_n]$ and 
$G$ a finite group acting linearly on $A$. Since
the construction of the finite minimal model of $A$ is functorial (it
is an application of the general construction in~\S\ref{ref-6.1-17})
 we obtain a finite minimal model for $\Lambda$ of the form
\begin{equation}
\label{ref-6.3-22}
\tilde{\Lambda}=\biggl(k\langle (x_S)_S\rangle\# G,d\biggr)
\end{equation}
(with $l=kG$).
The differential is $kG$-linear and on the variables $(x_S)_{S\neq \emptyset}$
it is given by the formula \eqref{ref-6.2-21}.
\subsection{Crossed products for cyclic groups}
\label{ref-6.4-23}
In this section we specialize to the case where $G$ is the cyclic
group $\ZZ/m\ZZ$. We now assume that $k$ is algebraically closed of
characteristic zero so that we may assume that $G$ acts diagonally on
$V=kx_1+\cdots+kx_n$.  Fix a primitive $m$'th root of unity~$\xi$. We
write $\chi_i$ for the character of $G$
\[
\chi_i(\bar{a})=\xi^{ai}
\]
and we let $e_i$ be the corresponding primitive idempotent
\[
e_i=\frac{1}{m}\sum_{a=0}^{m-1} \chi_i(\bar{a})\bar{a}=\frac{1}{m}\sum_{a=0}^{m-1}\xi^{ai} \bar{a}
\]
so that we have
\[
kG=\sum_{i=0}^{m-1} ke_i
\]
We assume that $G$ acts on the variable $x_i$  by the character $\chi_{a_i}$ for some $0\le a_i\le m-1$. Thus
in $\Lambda=A\# G$ with $A=k[x_1,\ldots,x_n]$ we have the relation
\[
\bar{a}\cdot x_i=\xi^{aa_i} x_i\cdot \bar{a}
\]
which implies
\begin{align*}
e_j\cdot x_i=x_i \cdot e_{j+a_i}
\end{align*}
(from now we tacitly reduce indices modulo $m$). Specializing \eqref{ref-6.3-22} we obtain a finite minimal model $\tilde{\Lambda}$ for $\Lambda$ which is 
freely generated over $l$ by the variables $(x_S)_S$ subject to the relations
\[
e_j\cdot x_S=x_S \cdot e_{j+d(S)}
\]
where 
\[
d(S)=\sum_{i\in S} a_i
\]
and with differential as in \eqref{ref-6.2-21}.
\subsection{McKay quiver description}
We reformulate the results of the previous section in quiver language.
Let $Q$ be the
quiver with vertices $\{0,\cdots,m-1\}$ (taken modulo $m$) and arrows
$x_{j,i,j+a_i}$ for $i=1,\ldots,n$, $j=0,\ldots,m-1$ starting at $j$
and ending at $j+a_i$.  Then $\Lambda$ is a quotient of the path algebra of
$Q$ where the arrow $x_{i,j,j+a_i}$ is sent to
$e_jx_i=x_ie_{j+a_i}=e_jx_ie_{j+a_i}\in \Lambda$. The relation $x_kx_l=x_lx_l$ in
$\Lambda$ leads to the relation
\begin{equation}
\label{ref-6.4-24}
x_{j,k,j+a_k}\cdot x_{j+a_k,l,j+a_k+a_l}=x_{j,l,j+a_l}\cdot x_{j+a_l,k,j+a_k+a_l}
\end{equation}
on $Q$. The path algebra of $Q$ module all such relations is precisely $\Lambda$. 

To describe the algebra $\tilde{\Lambda}$ we introduce the graded quiver
$\tilde{Q}$ with the same vertices as $Q$ and with arrows
$x_{j,S,j+d(S)}$ of homological degree $-|S|+1$ going from $j$ to
$j+d(S)$. The differential on $k\tilde{Q}$ is given by
\begin{equation}
\label{ref-6.5-25}
dx_{j,S,j+d(S)}=\sum_{S=A\coprod B,A\neq\emptyset,B\neq\emptyset}(-1)^{|A|-1} \epsilon_{A,B} x_{j,A,j+d(A)} \cdot x_{j+d(A),B,j+d(S)}
\end{equation}
Let $\tilde{Q}^0$ be the quiver obtained from $\tilde{Q}$ by dropping
all arrows adjacent to $0$. Thus the arrows in $\tilde{Q}^0$ are
of the form $x_{j,S,j+d(S)}$ with $j\neq 0$, $j+d(S)\neq 0$. Specializing
Theorem \ref{ref-5.1.1-5} to the current situation we obtain the
following result which we will formulate in the complete setting
for variety. Note that by Remark \ref{ref-6.1.2-20} this makes little difference.
\begin{propositions}
\label{ref-6.5.1-26}
Let $A=k[x_1,\ldots,x_n]$ and assume that the cyclic group $G=\ZZ/m\ZZ$ acts
linearly on $A$ with weights $(\xi^{a_1},\ldots,\xi^{a_n})$ satisfying the
additional properties $\sum_ia_i\cong 0\text{ mod } 0$ and $\gcd(a_i,m)=1$.
Then 
\[
\underline{\MCM}(\hat{A}^G)\cong \Perf(k\tilde{Q}^0,d)/\langle
(S_i)_{i=1,\ldots,m-1}\rangle
\]
where the differential is given by \eqref{ref-6.5-25} (taking into
account that arrows adjacent to the vertex $0$ should be suppressed
on the righthand side). The DG-algebra $(k\tilde{Q}^0,d)$ has finite
dimensional cohomology and 
\[
H^0(k\tilde{Q}^0,d)\cong \hat{\Lambda}/\hat{\Lambda} e_0\hat{\Lambda}
\]
\end{propositions}
\begin{proof} This is a straight translation of Theorem \ref{ref-5.1.1-5}.
The conditions on the numbers $(a_i)_i$ are to insure that $\hat{A}^G$ is
Gorenstein and has an isolated singularity. 
\end{proof}
\begin{examples} We discuss an example that occurred in \cite{IY,
KMVdB}. Assume that $G=\ZZ/2\ZZ$ acts with weights $(-1,-1,-1,-1)$
on $k[x_1,x_2,x_3,x_4]$. The resulting McKay quiver looks like
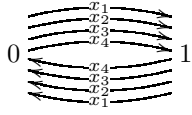
\[
\xymatrix{
0 
\ar@<3ex>@/^0.5em/[rr]|{x_1}
\ar@<2ex>@/^0.5em/[rr]|{x_2}
\ar@<1ex>@/^0.5em/[rr]|{x_3}\ar@/^0.5em/[rr]|{x_4}&&
1
\ar@<3ex>@/^0.5em/[ll]|{x_1}
\ar@<2ex>@/^0.5em/[ll]|{x_2}
\ar@<1ex>@/^0.5em/[ll]|{x_3}\ar@/^0.5em/[ll]|{x_4}&&
}
\]
with relations
\[
x_ix_j=x_jx_i
\]

The quiver $\tilde{Q}^0$ contains a single vertex $1$ and loops
$x_{ij}$, $i<j$ and $x_{1234}$ respectively of degree $-1$ and $-3$. The differential is given by
\[
dx_{1234}=-[x_{12},x_{34}]-[x_{23},x_{14}]+[x_{24},x_{13}]
\]
where as usual $[-,-]$ denotes the graded commutator. 
\end{examples}
\subsection{The minimal case}
In this section we impose the additional condition.
\begin{equation}
\label{ref-6.6-27}
\sum_i a_i=m 
\end{equation}
It was proved in \cite{AIR} that if condition \eqref{ref-6.6-27} 
holds then $\underline{\MCM}(\hat{A}^G)$ has a description as a higher cluster
category. In this section we indicate how to prove this starting from
Proposition \ref{ref-6.5.1-26}.

We order the vertices of $\tilde{Q}^0$ by their label. We we say that
an arrow $x_{j,S,j+d(S)}$ is ascending if $j+d(S)>j$ (recall that the
sum $j+d(S)$ is reduced $\mod m$). Otherwise we say that the arrow is
descending.

Let $P$ be the quiver with the same vertices as $\tilde{Q}^0$ but only
with ascending arrows of the form $x_{j,i,j+a_i}$ subject to relations
\[
x_{j,k,j+a_k}\cdot x_{j+a_k,l,j+a_k+a_l}=x_{j,l,j+a_l}\cdot x_{j+a_l,k,j+a_k+a_l}
\]
Let $C$ be the resulting path algebra. As in \cite{Keller11,VdB29} we
define the inverse dualizing complex of $C$ as
\[
\Phi_C=\RHom_{C^e}(C,C\otimes C)
\]
\begin{propositions}
\label{ref-6.6.1-28}
There is a quasi-isomorphism of DG-algebras
\begin{equation}
\label{tensor}
(k\tilde{Q}^0,d)\cong T_C(\Theta_{C}[n-1])
\end{equation}
In particular $
\underline{\MCM}(\hat{A}^G)
$
is a generalized $n-1$-cluster category.
\end{propositions}
The fact that $\underline{\MCM}(\hat{A}^G)$ has an $n-1$-cluster tilting object
has already been established by Iyama \cite[Thm 5.2.1]{iyama}. So we only need
to be concerned with proving \eqref{tensor}.
The rest of this section will be devoted to this.

Let $\tilde{P}$ be the quiver with the same vertices as $\tilde{Q}^0$
but only with ascending arrows.  Let $D$ be the $k\tilde{P}$ bimodule
generated by the descending arrows. Then we have
\[
k\tilde{Q}^0=T_{k\tilde{P}}D
\]
Furthermore if we look at the formula \eqref{ref-6.5-25} we see that both
$k\tilde{P}$ and $D$ are closed under the differential $d$ (it is here
that \eqref{ref-6.6-27} is used).

Hence we need to prove two things
\begin{enumerate}
\item The natural map $k\tilde{P}\r C$ is a quasi-isomorphism.
\item
$D\cong \Phi_{k\tilde{P}}[n-1]$
\end{enumerate}
The first statement follows from the following lemma
\begin{lemmas}
The algebra $C$ is Koszul over $l^0\overset{\text{def}}{=}l/le_0$. Furthermore $k\tilde{P}$ is the 
minimal model for $C$ constructed in \S\ref{ref-6.1-17}.
\end{lemmas}
\begin{proof}
We give $A$ an additional grading by putting $\deg x_i=a_i$. It is easy to
see that the minimal resolution of $S_i$ for $C$ is obtained by resolving $k$
as a graded $A$-module for the Adams grading 
and then truncating in degrees $\le m-i$ for the additional grading. It is clear
that this resolution is linear for the Adams grading. 

Given the explicit form of the relations of $C$ one verifies immediately that
$k\tilde{P}$ is indeed the finite minimal form constructed in \S\ref{ref-6.1-17}. 
\end{proof}
Now we prove the second statement.  We first discuss some generalities
on bimodules and differentials.

Let $E$ be a DG-algebra and assume that $P$,$Q$ are graded DG-$E$-bimodules. We write 
$
P\otimes_{E^w} Q
$
for the quotient of $P\otimes_E Q/[E,P\otimes_E Q]$. The dual of an $E$-bimodule $Q$
is by definition
$
Q^\ast=\Hom_{E^e}(Q,{}_EE\otimes E_E)
$.
This is again a DG-$E$-bimodule through the surviving inner bimodule structure
on $E\otimes E$. 

If $\omega\in (P\otimes_{E^e} Q)_n$ then we obtain an induced map of
degree $n$
\[
\omega^\circ:P^\ast\r Q:\phi\mapsto (-1)^{|\phi(\omega')'|(|\omega''|+|\phi(\omega')''|)}\phi(\omega')''\omega''\phi(\omega')'
\]
We say that $\omega$ is non-degenerate if $P$, $Q$ are finitely generated
projective $E$-bimodules and $\omega^\circ$ induces an isomorphism $P^\ast\cong Q[n]$. 
If $d\omega=0$
then $\omega^\circ$ is a morphism of DG-bimodules. 

\medskip

As $E$-bimodule $E$ is quasi-isomorphic to the cone of
\[
0\r \Omega_{E/l}\xrightarrow{\sigma} E\otimes_l E\r 0
\]
where $\Omega_{E/l}$ is the bimodule of non-commutative differentials. As
a graded $E$-bimodule $\Omega_{E/l}$ is generated by elements $Db$ subject to the standard
relations.\footnote{Using $db$ would lead to confusion with the differential
on $E$} One has $\sigma(Db)=b\otimes 1-1\otimes b$. The operator $D$
has degree zero, whence $d(Db)=D(db)$. 

We denote the cone of $\sigma$ by $\tilde{\Omega}_{E/l}$. If $b\in
E$ then $Db$ considered as an element of degree $|b|-1$ of
$\tilde{\Omega}_{E/l}$ is written as $\tilde{D}b$.

Hence $\tilde{\Omega}_{E/l}$ is generated by $\tilde{D}b$ and an
$l$-central element $g$ of degree zero (this represents the generator
$1\otimes 1$ of $E\otimes_l E$). The differential on
$\tilde{\Omega}_{E/l}$ is given by
\begin{equation}
\label{ref-6.7-29}
\begin{aligned}
dg&=0\\
d(\tilde{D}b)&=-\tilde{D}(db)+[b,g]
\end{aligned}
\end{equation}
We compute $\tilde{\Omega}_{k\tilde{P}/l^0}$. For $j<j+d(S)$ write
$\tilde{x}_{j,d(S),j+d(S)}=\tilde{D} x_{j,S,j+d(S)}$. Thus $|\tilde{x}_{j,d(S),j+d(S)}|=
-|S|$.
Furthermore we
introduce loops $\tilde{x}_{j,\emptyset,j}$ at $j$ of degree zero which correspond to
$-ge_j$. Then $\tilde{\Omega}_{k\tilde{P}/l^0}$ is the free bimodule
with generators $\tilde{x}_{j,d(S),j+d(S)}$ for $j\le j+d(S)$,
$S\subsetneq [n]$ with $[n]=\{1,\ldots,n\}$. A simple computation using
the formulas \eqref{ref-6.7-29}
yields
\begin{multline}
\label{ref-6.8-30}
d(\tilde{x}_{j,S,j+d(S)})=\sum_{S=A\coprod B,B\neq \emptyset} (-1)^{|A|}\epsilon_{A,B}  
\tilde{x}_{j,A,j+d(A)}\cdot
x_{j+d(A),B,j+d(S)}\\-\sum_{S=A\coprod B,A\neq \emptyset} \epsilon_{A,B} x_{j,A,j+d(A)}\cdot
\tilde{x}_{j+d(A),B,j+d(S)}
\end{multline}
We define an element $\omega\in \tilde{\Omega}_{k\tilde{P}/l^0}\otimes_{k\tilde{P}^e}D$ by
the following formula
\[
\omega=\sum_{S\subsetneq [n]),j\le j+d(S)}
(-1)^{|S|-1}\epsilon_{S,S^c}\tilde{x}_{j,S,j+d(S)}\otimes x_{j+d(S),S^c,j}
\]
One verifies
\begin{enumerate}
\item $|\omega|=-n+1$.
\item $d\omega=0$.
\item $\omega$ is non-degenerate.
\end{enumerate}
It follows that $\omega$ defines an isomorphism 
$\omega^\circ:\Theta_{k\tilde{P}}=\tilde{\Omega}_{k\tilde{P}/l^0}^\ast\r D[-n+1]$.
Hence the proof is complete.
\def\cprime{$'$} \def\cprime{$'$} \def\cprime{$'$} \def\cprime{$'$}
\ifx\undefined\bysame
\newcommand{\bysame}{\leavevmode\hbox to3em{\hrulefill}\,}
\fi

\end{document}